\documentclass[11pt]{amsart}
\usepackage{amsmath,amssymb,amsfonts,url}
\numberwithin{equation}{section}

\newtheorem{thm}{Theorem}[section]
\newtheorem{lem}{Lemma}[section]

\newtheorem{prop}{Proposition}[section]

\begin{document}
\title[Monotonicity Theorem]{A local parabolic monotonicity formula on Riemannian manifolds}

\author{Eduardo V. Teixeira}
\address{Universidade Federal do Cear\'a\\ Departamento de
Matem\'atica\\ Av. Humberto Monte, s/n, Campus do Pici - Bloco 914
\\
Fortaleza-CE, Brazil.  CEP 60.455-760} \email{eteixeira@pq.cnpq.br}

\author{Lei Zhang}
\address{Department of Mathematics\\
        University of Alabama at Birmingham\\
        1300 University Blvd, 452 Campbell Hall\\
        Birmingham, Alabama 35294-1170}
\email{leizhang@math.uab.edu}
\thanks{E. Teixeira acknowledges support from NSF and CNPq.
L. Zhang is supported in part by NSF Grant 0600275 (0810902)}

\date{\today}

%%%%%%%%%%%%%%%%%%%%%%%%%%%%%%%%%%%%%%%%%%%%%
\begin{abstract}
In this article we establish a local parabolic almost monotonicity
formula for two phase free boundary problems on Riemannian
manifolds, which is an extension of a work of Edquist-Petrosyan.
\end{abstract}
%%%%%%%%%%%%%%%%%%%%%%%%%%%%%%%%%%%%%%%%%%%%%

\maketitle

\section{Introduction}

In the theory of two-phase free boundary problems, it is well
established that regularity of the interface is closely related to
asymptotic behavior of solution near the free boundary. In 1984
Alt, Caffarelli and Friedman \cite{alt} established a monotonicity
formula to describe the interaction of the two pieces of the
solution on each side of the free boundary. This formula has been
extremely powerful in the regularity theory and it reads as
follows: Let $u_1,u_2$ be two non-negative continuous functions in
$B_1$ (the unit ball in $\mathbb R^n$) such that $\Delta u_i\ge 0$
($i=1,2$) are satisfied in distribution. Suppose $u_1\cdot u_2=0$
and $u_1(0)=u_2(0)=0$, then
$$\phi(r)=\frac{1}{r^4}
\int_{B_r}\frac{|\nabla u_1|^2}{|x|^{n-2}}dx \int_{B_r}\frac{|\nabla
u_2|^2}{|x|^{n-2}}dx$$ is monotone non-decreasing for $0<r<1$.

\bigskip

There have been different extensions of the theorem of
Alt-Caffarelli-Friedman under different contexts. For example,
Caffarelli \cite{caf1} established a monotonicity formula for
variable coefficient operators, Friedman-Liu \cite{friedman} have an
extension for eigenvalue problems. Another important extension has
been achieved by Caffarelli-Jerison-Kenig \cite{cjk} who replace
$\Delta u_i\ge 0$ by $\Delta u_i\ge -1$ ($i=1,2$). Under this new
assumption they prove that $\phi(r)$ is uniformly bounded for
$0<r<\frac 12$. This is called an "almost monotonicity formula".

Even though there is no monotonicity in the Caffarelli-Jerison-Kenig
formula, it does provide a control of $|\nabla u_1(0)|\cdot |\nabla
u_2(0)|$ if $u_1,u_2$ are both smooth at $0$. For many free boundary
problems, the control of $|\nabla u_1(0)|\cdot |\nabla u_2(0)|$
usually leads to important regularity results. Moreover, for some
real life problems such as the Prandtl-Batchelor problem
\cite{acker,batchlor1,batchlor2}) and some classical problems ( e.g.
see Shahgholian \cite{shahgholian}), the equations may be
inhomogeneous and we may not have $\Delta u_i\ge 0$ ($i=1,2$) on
each side of the free boundary. The ``almost monotonicity formula"
of Caffarelli-Jerison-Kenig is particularly useful in these
situations and has provided a theoretical basis for the regularity
theory of many new problems (see for example
\cite{cjk,shahgholian}).

For two-phase parabolic free boundary problems, Caffarelli
\cite{caffarelli5} established a monotonicity formula for two
sub-caloric functions: Let $u_1,u_2$ satisfy $\Delta
u_i-\partial_tu_i\ge 0$ in $\mathbb R^n\times (-1,0)$ ($i=1,2$),
$u_1u_2\equiv 0$, $u_1(0,0)=u_2(0,0)=0$. Let
\begin{equation}
\label{gxt} G_0(x,t)=\frac{1}{(4\pi t)^{\frac
n2}}e^{-\frac{|x|^2}{4t}}
\end{equation}
be the fundamental solution of the heat equation in $\mathbb R^n$,
then
\begin{equation}\label{monotone}
\phi(r)=\frac{1}{r^4}\int_{\mathbb R^n}\int_{-r^2}^0|\nabla
u_1|^2G_0(x,-s)dsdx\int_{\mathbb R^n}\int_{-r^2}^0|\nabla
u_2|^2G_0(x,-s)dsdx
\end{equation}
is monotone non-decreasing provided that $u_1,u_2$ have reasonable
growth at infinity.

\par

Clearly this monotonicity formula for parabolic free boundary
problems is in correspondence with the Alt-Caffarelli-Friedman
formula. Later, Edquist-Petrosyan \cite{edquist} derived the "almost
monotonicity formula" for $u_1,u_2$ satisfying $\Delta
u_i-\partial_tu_i\ge -1$ instead of being sub-caloric. Similar to
the Caffarelli-Jerison-Kenig formula, Edquist-Petrosyan proved the
bound of $\phi(r)$ for $0<r<\frac 12$.

\par
A common feature of all monotonicity and almost monotonicity
formulas aforementioned is that they are all designed for problems
within Euclidean spaces. From theoretical and application viewpoints
it is natural to consider some free boundary problems on Riemannian
manifolds. Indeed, it has been pointed out by Caffarelli and Salsa
in \cite{salsa1} that the tools developed for free boundary problems
on Euclidean spaces should have their counterparts for free boundary
problems on manifolds (page ix of the introduction). The analogs of
Alt-Caffarelli-Friedman monotonicity formula and
Caffarelli-Jerison-Kenig almost monotonicity formula have been
developed for the Laplace-Beltrami operator by the authors in
\cite{tz1}. The purpose of this article is to derive the analogue of
Edquist-Petrosyan formula on Riemannian manifolds.  In forthcoming
works we shall use these formulas to discuss free boundary problems
on Riemannian manifolds.

\par

Let $(M,g)$ be a Riemannian manifold of dimension $n\ge 2$, let
$B(p,\delta_p)$ be a geodesic ball around $p$ with radius
$\delta_p=\min\{1,\mbox{inj}_p\}$ ($\mbox{inj}_p$ is the
injectivity radius at $p$). We shall use the following cut-off
function $\chi$ supported in $B(p, \delta_p)$:
$$\chi\equiv 1, \mbox{ in } B(p, \delta_p/4),\quad \chi\equiv 0, \mbox{ in } B(p,
\delta_p/2). $$ Let $R_m$ denote the curvature tensor. Assume
\begin{equation}\label{feb10e1}
|R_m|+|\nabla_gR_m|\le \Lambda.
\end{equation}
In this article, if we do not mention the dependence of a given
constant, it is implied that this constant is either universal or
depends only on $n$, $\Lambda$ and $\chi$.

 Let $Q_r^-(p)=B_p(r)\times
(-r^2,0)$ for $0<r<\delta_p/2$. If $u_+,u_-\in
H^1_{loc}(Q_{\delta_p}^-(p))$, we define $w_{\pm}$ as
$w_{\pm}=u_{\pm}\chi$ and let
$$\phi(r)=\frac
1{r^4}\int\!\!\!\int_{S_r(p)}|\nabla_gw_+|^2G(x,-s)dV_gds\int\!\!\!\int_{S_r(p)}|\nabla_gw_-|^2G(x,-s)dV_gds$$
where $G(x,t)=\frac 1{(4\pi t)^{\frac n2}}e^{-\frac{d(x,p)^2}{4t}}$,
$S_r=B(p,\delta_p)\times (-r^2,0)$, $dV_g=\sqrt{\mbox{det}(g)}dx$.
Our main result is:
\begin{thm} \label{thm1} Let $u_{\pm}$ be non-negative continuous
functions in $Q_{\delta_p}^-(p)$ that satisfy $u_+\cdot u_-=0$ and
$$
(\Delta_g-\partial_t)u_{\pm}\ge -1,\quad \mbox{in}\quad Q_p^-(1)$$
in the weak sense. Then $u_{\pm}\in H^1_{loc}(Q_{\delta_p}^-(p))$
and there exists $C>0$ such that
$$\phi(r)\le C(1+\int\!\!\!\int_{Q_{\delta_p}^-}u_+^2dV_gds+\int\!\!\!\int_{Q_{\delta_p}^-}
u_-^2dV_gds)^2 $$ for all $r\in (0, \frac 12\delta_p)$.
\end{thm}

The proof of Theorem \ref{thm1} is along the lines of \cite{cjk} and
\cite{edquist}. However, since equations, integrals and kernels are
defined on a Riemannian manifolds, many perturbation terms have to
be properly controlled in different situations. For example, we need
to derive a ``perturbed" version of Beckner-Kenig-Pipher inequality
in Proposition \ref{09428p1}. As the reader will see, it is crucial
in our analysis to have a quantitative estimate of the perturbation
of the eigenvalues in this key inequality. In accordance to
\cite{cjk} and \cite{edquist}, assuming Holder continuity of
solution, we can infer a more precise control of the functional
$\phi$. This is the content of Theorem \ref{thm2} we present at the
end of the paper. With our ``perturbed" version of
Beckner-Kenig-Pipher inequality in hand, the proof of Theorem
\ref{thm2} becomes very similar to the corresponding theorem in
\cite{edquist} and is therefore omitted.

\section{The proof of the monotonicity formula}

Let $U(x,t)$ be the heat kernel in the neighborhood of $p$, then for
$t$ small we have
$$U(x,t)=(4\pi t)^{-\frac n2}e^{-\frac{d^2(x)}{4t}}\bigg
(\sum_{i=0}^{\infty}\phi_i(x)t^i\bigg ) $$ (see \cite{schoenyau}
P109) where $\phi_i$ are smooth functions of $x$, and
$$\phi_0=\mbox{det}(g)^{-\frac 14}=1+O(d(x)^2). $$
Since we consider a neighborhood of the origin, $G$ and $U$ are
clearly comparable in this neighborhood. As a consequence, $G$ can
be replaced by $U$ in Theorem \ref{thm1}. In the following we will
mainly use $U$ in our proof.

Let
$$A_k^{\pm}=\int\!\!\!\int_{S_{4^{-k}}}|\nabla_gw_{\pm}|^2U(x,-s)dV_gds,
\quad b_k^{\pm}=4^{4k}A_k^{\pm}.$$

To prove Theorem \ref{thm1} it is enough to prove the bound of
$\phi(r)$ for $0<r<\delta$ where $\delta$ is a small constant
depending on $n,\Lambda$ and $\chi$. The bound for $r>\delta$ is
obvious. Therefore in the proof we only focus on the estimates of
$A_k^{\pm}$ for $k$ sufficiently large.

We shall prove the following two key propositions:
\begin{prop}\label{keyprop1}
There exist $C_0,C_1>0$ such that for $k\ge k_0(n,\Lambda,\chi)$, if
$b_k^{\pm}\ge C_0$,
$$4^4A_{k+1}^+A_{k+1}^-\le
A_k^+A_k^-(1+\delta_k)$$ where
$\delta_k=C_1(\frac{1}{\sqrt{b_k^+}}+\frac{1}{\sqrt{b_k^-}}+4^{-2k})$.
\end{prop}

\begin{prop}\label{keyprop2}
There exists $\epsilon\in (0,1)$ such that for $k\ge
k_0(n,\Lambda)$, if $b_k^{\pm}\ge C_0$ and $4^4A_{k+1}^+>A_k^+$ then
$A_{k+1}^-\le (1-\epsilon)A_k^-$.
\end{prop}

Theorem \ref{thm1} follows from Proposition \ref{keyprop1} and
Proposition \ref{keyprop2} by standard argument in \cite{cjk,tz1}.
Here we note that we shall always assume $u_{\pm}$ to be smooth, as
$u_{\pm}$ can be mollified to $u^{\epsilon}_{\pm}$ such that
necessary inequalities for $u^{\epsilon}_{\pm}$ can be obtained
first. Therefore, the conclusion for $u_{\pm}$ can be obtained by
passing $u^{\epsilon}_{\pm}$ to the limit. This part of the argument
is standard and is omitted. The interested readers may look into
\cite{tz1} for reference.

The following estimate is important for the proof of both
Proposition \ref{keyprop1} and Proposition \ref{keyprop2}:
\begin{eqnarray}
&&\int\!\!\!\int_{S_r}|\nabla_gw_{\pm}|^2U(x,-s)dV_gds \nonumber \\
&\le &C_Mr^4+C_Mr^2(\int_{\mathbb
R^n}w_{\pm}^2(x,-r^2)U(x,r^2)dV_g)^{\frac 12}\nonumber \\
\label{322e1}
 &&+\frac 12\int_{\mathbb R^n}w_{\pm}^2(x,-r^2)U(x,r^2)dV_g.
\end{eqnarray}
where $C_M$ depends on $\int\!\!\!\int_{S_1}(u_+^2+u_-^2)dV_gds$.

\medskip

\noindent{\bf Proof of (\ref{322e1}):} From
$(\Delta_g-\partial_s)u_{\pm}\ge -1$, we have, by standard
computation
$$
(\Delta_g-\partial_s)(w_{\pm}^2) \ge -2w_{\pm}\chi
-4|\nabla_gu_{\pm}|\cdot |\nabla_g\chi
|w_{\pm}-2u_{\pm}|\Delta_g\chi |w_{\pm}+2|\nabla_gw_{\pm}|^2.
$$
From the above we have
\begin{eqnarray*}
&&2\int\!\!\!\int_{S_r}|\nabla_gw_{\pm}|^2d\nu\le
\int\!\!\!\int_{S_r}(\Delta_g-\partial_s)(w_{\pm}^2)d\nu\\
&&+2\int\!\!\!\int_{S_r}w_{\pm}d\nu+4\int\!\!\!\int_{S_r}(|\nabla_gu_{\pm}|
|\nabla_g\chi |+u_{\pm}|\Delta_g \chi |)w_{\pm}d\nu\\
&=&I_1+I_2+I_3.
\end{eqnarray*}
Here we use the notation: $d\nu=U(x,-s)dV_gds$. We shall also use
$d\nu^s=U(x,-s)dV_g$.

Using $w_{\pm}\ge 0$ and integration by parts we obtain
\begin{equation}\label{i1}
I_1\le \int_{\mathbb R^n}w_{\pm}^2(x,-r^2)d\nu^{-r^2}.
\end{equation}
 Note that we used
$$\lim_{\epsilon\to
0^+}\int_{\mathbb
R^n}w_{\pm}^2(x,-\epsilon)U(x,\epsilon)dV_g=-w_{\pm}^2(0,0)\le 0.
$$

 To estimate $I_2$ we use the following equation
easy to be verified by direct computation:
$$(\Delta_g-\partial_s)(w_{\pm}-s)\ge
2\nabla_gu_{\pm}\nabla_g\chi+u_{\pm}\Delta_g\chi. $$ For $s_2\le
s_1\le 0$ we integrate the above to obtain
\begin{eqnarray}
&&\int_{\mathbb R^n}(w_{\pm}-s_1)d\nu^{s_1} \le \int_{\mathbb
R^n}(w_{\pm}-s_2)d\nu^{s_2}\nonumber \\
\label{323e1} &&+\int_{s_2}^{s_1}\bigg (\int_{\mathbb
R^n}(2|\nabla_gu_{\pm}| |\nabla_g\chi |+u_{\pm}|\Delta_g\chi
|)d\nu^s \bigg )ds.
\end{eqnarray}

For $0\ge s_1\ge -r^2\ge s_2\ge -4r^2$, since the support of $\nabla
\chi$ or $\Delta \chi$ stays away from the origin, it is elementary
to obtain
\begin{eqnarray}
&&\int_{\mathbb R^n}w_{\pm}(x,s_1)d\nu^{s_1}\le \inf_{s\in
[-4r^2,-r^2]}\int_{\mathbb R^n}w_{\pm}(x,s)d\nu^s \nonumber \\
\label{323e2} && +C_1(\Lambda,n)r^2+C_2(M,N)r^N.
\end{eqnarray}
where $M$ is the $L^2$ norm of $u_{\pm}$ on $Q_1^-$, $N$ is a large
number. Therefore $I_2$ satisfies
\begin{equation}\label{i21}
I_2\le 2r^2\inf_{s\in [-4r^2,-r^2]}\int_{\mathbb
R^n}w_{\pm}(x,s)d\nu^s+C(N,M)r^4.
\end{equation}
If we further use Cauchy's inequality we have
$$I_2\le C(M,\Lambda)r^4+\inf_{s\in
[-4r^2,-r^2]}\int_{\mathbb R^n}w_{\pm}^2(x,s)d\nu^{s}.
$$
 The estimate of $I_3$ is similar.
 \begin{equation}\label{i3}
I_3\le
C(N)r^N\int\!\!\!\int_{Q_{\frac{\delta_p}2}^-}(|\nabla_gu_{\pm}|^2+|u_{\pm}|^2)\le
C(N,M)r^N\int\!\!\!\int_{Q_{\frac{\delta_p}2}^-}u_{\pm}^2.
\end{equation}
(\ref{322e1}) can be obtained easily from (\ref{i1}), (\ref{i21})
and (\ref{i3}).
 As a consequence, the following
estimates also hold:
\begin{equation}\label{324e1}
\int\!\!\!\int_{S_r}|\nabla_gw_{\pm}|^2d\nu\le C_Mr^4+C_M\inf_{s\in
[-4r^2,-r^2]}\int_{\mathbb R^n}w_{\pm}^2(x,s)d\nu^s.
\end{equation}

\begin{equation}\label{324e2}
\int\!\!\!\int_{S_r}|\nabla_gw_{\pm}|^2d\nu\le
C_Mr^4+\frac{C_M}{r^2}\int\!\!\!\int_{S_{2r}\setminus
S_r}w_{\pm}^2(x,s)d\nu.
\end{equation}

In the following, we shall always re-scale $w_{\pm}$ as follows:
$$\bar w_{\pm}(y,s)=\frac 1{r^2}w_{\pm}(ry,r^2s).$$
$\bar u_{\pm}$ is understood similarly. Correspondingly we let
$$\bar
g_{ij}(\cdot )=g_{ij}(r\cdot)=\delta_{ij}+O(r^2)|\cdot |^2$$ be the
re-scaled metric, $d\bar \nu$ and $d\bar \nu^s$ are defined as
$$d\bar \nu=r^nU(ry,-r^2s)dV_{\bar g}ds,\quad
d\bar \nu^s=r^nU(ry,-r^2s)dV_{\bar g},$$ then $\bar u_{\pm}$ satisfy
$$(\Delta_{\bar g}-\partial_s)\bar u_{\pm}\ge -1.$$
We use $\Omega_{\pm}$ to represent the set where $w_{\pm}$ is
positive. The corresponding set for $\bar w_{\pm}$ is $\bar
\Omega_{\pm}$.

\begin{lem}\label{poincare} For $r<\delta_p$ and $s=-\frac 12$
\begin{equation}\label{53e2}
\log \frac{1+O(r^2)}{|\bar w_{\pm}|_{\bar \nu^s}}\int_{\mathbb
R^n}\bar w_{\pm}^2d\bar\nu^s\le 2(1+O(r^2))\int_{\mathbb
R^n}|\nabla_{\bar g}\bar w_{\pm}|^2d\bar\nu^s.
\end{equation}
\end{lem}

\noindent{\bf Proof of Lemma \ref{poincare}:}  We use $d\nu_0$ to
represent the Gauss measure in Euclidean spaces:
$$d\nu_0=\frac{1}{(2\pi)^{\frac{n}2}}e^{-\frac{|x|^2}{2}}dx.$$
We perform two transformations on the three integral terms in
(\ref{53e2}):
$$
 |\bar w_{\pm}|_{\bar \nu^s},\quad \int_{\mathbb
R^n}\bar w_{\pm}^2d\bar\nu^s,\quad \int_{\mathbb R^n}|\nabla_{\bar
g}\bar w_{\pm}|^2d\bar\nu^s$$
 to reduce them to the Euclidean case.
First, using $y=\phi_1(x)$ ($\phi_1$ to be determined) we have
\begin{eqnarray}
&&|\nabla_{\bar g}\bar w_{\pm}(x)|^2 =\bar g^{ij}(x)\frac{\partial
\bar w_{\pm}}{\partial x_i}\frac{\partial \bar w_{\pm}}{\partial
x_j}\nonumber\\
\label{phi1} &=&\bar g^{ij}(\phi_1^{-1}(y))\frac{\partial
y^m}{\partial x_i} \frac{\partial y^l}{\partial x_j}\frac{\partial
\bar w_{\pm}}{\partial y_m}\frac{\partial \bar w_{\pm}}{\partial
y_l}
\end{eqnarray}
Here repeated indices imply summation. Since $\bar g_{ij}(x)$ is
symmetric and $\bar g_{ij}(x)=\delta_{ij}+O(r^2|x|^2)$ for $|x|\le
r^{-1}\delta_2$ ( $\delta_2$ small), we can choose $\phi_1$ so that
$$\frac{dy}{dx}=\bigg (\bar g_{ij}(x)\bigg )^{\frac 12}=\bigg (\delta_{ij}+O(r^2|x|^2)\bigg ).$$
Consequently
$$\bar g^{ij}(\phi_1^{-1}(y))\frac{\partial y^m}{\partial x_i}
\frac{\partial y^l}{\partial x_j}=\delta^{ml} $$ and (\ref{phi1})
becomes $|\nabla_{\bar g}\bar
w_{\pm}(x)|^2=\sum_{i=1}^n|\frac{\partial \bar w_{\pm}}{\partial
y_i}|^2$. Moreover the Jacobian of the mapping is of the order
$1+O(r^2|y|^2)$ for $|y|\le \delta_3r$ with $\delta_3$ small. With
this $\phi_1$ we combine the Jacobian with the heat kernel:
$$d\bar \nu_y:=d\bar \nu J_{\phi_1}=
r^nU(r\phi^{-1}(y),-r^2s)J_{\phi_1}dV_{\bar g}ds.$$ Using the
definition of $\phi_1$ we now have (recall that $s=-\frac 12$)
\begin{equation}\label{953e1}
r^nU(r\phi^{-1}(y),-r^2s)J_{\phi_1}dV_{\bar g}=\frac 1{(2\pi)^{\frac
n2}}e^{-\frac{|y|^2}2}(1+O(r^2)|y|^2)dy.
\end{equation}
We use $d\bar \nu_y^s$ to denote $r^nU(rx,-r^2s)J_{\phi_1}dV_{\bar
g}$. With these notations, the integral forms in (\ref{53e2}) become
\begin{eqnarray}
|\bar w_{\pm}|_{\bar \nu^s}=\int_{\mathbb R^n}\bar
w_{\pm}(\phi_1^{-1}(y))d\bar \nu_y^s\nonumber\\
\int_{\mathbb R^n}\bar w_{\pm}^2(x)d\bar \nu^s= \int_{\mathbb
R^2}\bar w_{\pm}^2(\phi_1^{-1}(y))d\bar \nu_y^s\nonumber\\
\label{54e4} \int_{\mathbb R^n}|\nabla_{\bar g}\bar w_{\pm}|^2d\bar
\nu^s= \int_{\mathbb R^n}|\nabla \bar
w_{\pm}(\phi_1^{-1}(y))|^2d\bar \nu_y^s.
\end{eqnarray}

The purpose of the second transformation is to make $d\bar \nu^s_y$
as close to the Gauss measure on Euclidean spaces as possible. To
this end we write $d\bar \nu_y^s$ as
$$
d\bar \nu_y^s=\frac{1}{(2\pi)^{\frac
n2}}e^{-\frac{|y|^2}2}(1+A(y))dy$$ where
$$|A(y)|\le Cr^2(1+|y|)^2,\quad |DA(y)|\le Cr^2(1+|y|),\quad |y|\le
\delta_3r^{-1}.$$ where $\delta_3$ is a small number. The second
transformation is defined as follows:
\begin{equation}\label{54e1}
y=z+\psi(z)
\end{equation}
 where $\psi$ satisfies $z\cdot \psi(z)=\ln(1+A(z))$ for $|z|>1$. It is easy to obtain from the estimate of $A$
that
\begin{equation}\label{psiz}
|\psi(z)|\le Cr^2|z|,\quad |D\psi(z)|\le Cr^2,\quad
1<|z|<\delta_4r^{-1}
\end{equation}
where $\delta_4$ is a small positive number. Then extend the
definition of $\psi$ to $B_1$ in such a way that both $|\psi|$ and
$|D\psi |$ are of the order $O(r^2)$ in $B_1$.

Using (\ref{54e1}) and (\ref{psiz}) we verify by direct computation
that
\begin{equation}\label{54e2}
\frac{1}{(2\pi)^{\frac n2}}e^{-\frac
12|y|^2}(1+A(y))J_{\frac{dy}{dz}}= \frac{1}{(2\pi)^{\frac
n2}}e^{-\frac 12|z|^2}(1+O(r^2)).
\end{equation}

Let $f_{\pm}(z)=\bar w_{\pm}(\phi^{-1}(y(z)))$, since the Jacobian
$J_{\frac{dy}{dz}}=1+O(r^2)$, the three integral terms in
(\ref{53e2}) are of the form (see (\ref{54e4})):
\begin{eqnarray}
|\bar w_{\pm}|_{\bar \nu^s}=\int_{\mathbb R^n}f_{\pm}(z)d\nu_0^s(1+O(r^2))\nonumber\\
\int_{\mathbb R^n}\bar w_{\pm}^2(x)d\bar \nu^s=\int_{\mathbb R^n}f_{\pm}^2(z)d\nu_0^s(1+O(r^2))\nonumber\\
\label{55e4} \int_{\mathbb R^n}|\nabla_{\bar g}\bar w_{\pm}|d\bar
\nu^s= \int_{\mathbb R^n}|\nabla f_{\pm}(z)|^2d\nu_0^s(1+O(r^2)).
\end{eqnarray}
Note that in the last equality, we used $\frac{dy}{dz}=id+O(r^2)$
where $id$ is the identity matrix. For $f_{\pm}$ we use the
Poincare's inequality on Euclidean spaces (see \cite{edquist}):
$$\log \frac{1}{|f_{\pm}|_{d\nu_0^s}}\int_{\mathbb
R^n}f_{\pm}^2d\nu_0^s\le 2\int_{\mathbb R^n}|\nabla
f_{\pm}|^2d\nu_0^s. $$
 Lemma \ref{poincare} follows from the equation above and (\ref{55e4}).
$\Box$

\bigskip

 The following two lemmas have analogues in
\cite{cjk,edquist} and their proofs are similar to their
counterparts in \cite{edquist}, we include the proofs here for the
convenience of the readers.
\begin{lem} \label{harp1}
 Let $w$ be $w_+$ or $w_-$, suppose
 $$\int\!\!\!\int_{\Omega\cap S_r}|\nabla_gw|^2d\nu=\alpha
r^4<\infty,\quad \int\!\!\!\int_{\Omega\cap S_{\frac
r4}}|\nabla_gw|^2d\nu\ge \frac{\alpha r^4}{256}.$$ Then there exists
$C>0$ and $C_1(n,\Lambda, M)>0$ such that if $\alpha>C$, $|\Omega
\cap S_{\frac r2}\setminus S_{\frac r4}|_{d\nu}\ge C_1r^2>0$, where
$\Omega$ is the set on which $w$ is positive.
\end{lem}

\noindent{\bf Proof of Lemma \ref{harp1}:}  Let $\bar w$ be $\bar
w_+$ or $\bar w_-$. Let $\bar \Omega$ be the set on which $\bar w$
is positive.  Then the assumptions become
$$\int\!\!\!\int_{\bar \Omega\cap S_1}|\nabla_{\bar g}\bar w|^2d\bar\nu
=\alpha,\quad \int\!\!\!\int_{\bar \Omega\cap S_{\frac
14}}|\nabla_{\bar g}\bar w|^2d\bar\nu\ge \frac{\alpha}{256}.
$$ We want
to show that if $\alpha$ is large, $|\bar \Omega \cap S_{\frac
12}\setminus S_{\frac 14}|_{d\bar\nu}>c_1(n)$. As a result of
(\ref{324e1}) we have
$$\frac{\alpha}{256}\le \int\!\!\!\int_{\bar \Omega\cap S_{1/4}}
|\nabla_{\bar g}\bar w|^2d\bar\nu\le C +C\inf_{s\in [-\frac
14,-\frac 1{16}]}\int_{\mathbb R^n}\bar w^2(x,s)d\bar\nu^s.$$

Therefore for $\alpha$ large
\begin{equation}\label{jun5e1}
\inf_{s\in [-\frac 14,-\frac 1{16}]}\int_{\mathbb R^n}\bar
w^2(x,s)d\bar\nu^s\ge \frac{\alpha}{512C}.
\end{equation}

From $\int\!\!\!\int_{S^1}|\nabla_{\bar g}\bar w|^2d\bar
\nu=\alpha$, we see that
$$\int_{\mathbb R^n}|\nabla_{\bar g}\bar w|^2d\bar \nu^s\le 16\alpha$$
except on a set of line measure no more than $\frac 1{16}$. So there
exists a set $E$ of line measure at least $\frac 18$ on $[-\frac
14,-\frac 1{16}]$ such that
$$\int_{\mathbb R^n}|\nabla_{\bar g}\bar w|^2d\bar \nu^s\le 16\alpha, \quad s\in E.
$$
For each $s\in E$, either $|\bar w|_{d\bar \nu^s}>\frac 12$, or
$|\bar w|_{d\bar \nu^s}\le \frac 12$. In this latter case we apply
Lemma \ref{poincare} and (\ref{jun5e1}) to get $|\bar w|_{d\bar
\nu^s}\ge c(n)$. Recall that we always assume $r$ to be small. So in
either case there exists $c(n)>0$ such that $|\bar w|_{d\bar
\nu^s}>c(n)$ for all $s\in E$. Therefore Lemma \ref{harp1} is
established by scaling. $\Box$

\begin{lem}\label{harp2} Let $\bar w$, $\bar \Omega$ be the same as
those in Lemma \ref{harp1}, assume $\mu\in (0,1)$ such that
$$|\bar \Omega\cap S_{\frac 12}\setminus S_{\frac 14}|_{d\bar \nu}\le (1-\mu)|S_{\frac
12}\setminus S_{\frac 14}|_{\bar \nu}.$$ Then there exists
$\lambda(\mu)\in (0,1)$ such that
$$\int\!\!\!\int_{S_{\frac 14}}|\nabla_{\bar g}\bar w|^2d\bar \nu\le
\lambda\int\!\!\!\int_{S_1}|\nabla_{\bar g}\bar w|^2d\bar \nu. $$
\end{lem}

\noindent{\bf Proof of Lemma \ref{harp2}:}
$$|\bar \Omega\cap S_{\frac 12}\setminus S_{\frac 14}|_{d\bar \nu}=\int_{-(\frac
14)^2}^{-(\frac 12)^2}|F(s)|ds$$ where $|F(s)|$ is the measure of
the positive set with respect to the $d\bar \nu^{s}$. We know that
$|F(s)|\le 1-\frac{\mu}2$ in a set $E\subset [-(\frac 12)^2,-(\frac
14)^2]$ with the line measure greater than or equal to
$\frac{\mu}2|S_{\frac 12}\setminus S_{\frac 14}|_{d\bar \nu}$.

Using (\ref{53e2}) we have, for small $r$ and $s\in E$, that
$$\int_{\mathbb R^n}\bar w^2(\cdot,s)^2d\bar \nu^s\le C\int_{\mathbb
R^n}|\nabla_{\bar g} \bar w(\cdot, s)|^2d\bar \nu^s$$

 If $\int\!\!\!\int_{\bar \Omega\cap S_{\frac 14}}|\nabla_{\bar
g}\bar w|^2d\nu\le \frac{\alpha}2$, there is nothing to prove.
Suppose this is not the case, then by using (\ref{324e1}) and the
largeness of $\alpha$, we have
$$\frac{\alpha}{4C}\le \int_{\mathbb R^n}\bar w^2(\cdot,s)d\bar \nu\quad \forall s\in [-\frac 14,-\frac{1}{16}].
$$
Specifically for $s\in E$ we have
$$\frac{\alpha}{4C}\le \int_{\mathbb R^n}\bar w^2(\cdot, s)d\bar \nu\le
C\int_{\mathbb R^n}|\nabla_{\bar g}\bar w(\cdot, s)|^2d\bar
\nu^s,\quad \forall s\in E.
$$

This implies that
$$\int_{\bar \Omega\cap S_{\frac 12}\setminus S_{\frac
14}}|\nabla_{\bar g}\bar w|^2d\bar \nu\ge \frac{\alpha |E|}{4C^2}.
$$ Lemma \ref{harp2} follows easily from the above.  $\Box$

\bigskip

The following proposition makes use of the two transformations used
in the proof of Lemma \ref{poincare}.

\begin{prop} \label{09428p1} Let $\bar \Omega_+^1\subset \mathbb R^n$ be the set where $\bar w_+(\cdot, -1)$ is positive.
$\bar \Omega_-^1$ is understood similarly. There exists $C>0$ such
that
$$\frac{\int_{\bar \Omega_+^1}|\nabla_{\bar g}\bar
w_+(\cdot,-1)|^2d\bar\nu^{-1}}{\int_{\bar \Omega_+}\bar
w_+(\cdot,-1)^2d\bar \nu^{-1}}+ \frac{\int_{\bar
\Omega_-^1}|\nabla_{\bar g}\bar w_-(\cdot,-1)|^2d\bar
\nu^{-1}}{\int_{\bar \Omega_-^1}\bar w_-(\cdot,-1)^2d\bar \nu^{-1}}
\ge 1-Cr^2.$$
\end{prop}

\noindent{\bf Proof of Proposition \ref{09428p1}:}

We make the two transformations as used in the proof of Lemma
\ref{poincare}. After the transformations, $\bar w_{\pm}(\cdot,-1)$
become $\tilde w_{\pm}$, $\bar \Omega_{\pm}^1$ become $\tilde
\Omega_{\pm}^1$. We still have $\tilde \Omega_+^1\cap \tilde
\Omega_-^1=\emptyset$ and $\tilde \Omega_+^1\cup
\tilde\Omega_-^1=B(0,\delta r^{-1})$ for some $\delta>0$ small.
$\tilde w_{\pm}$ are supported in $\tilde \Omega_{\pm}^1$,
respectively. Moreover, by the same estimates as in the proof of
Lemma \ref{poincare} we have
\begin{equation}\label{55e1}
 \int_{\bar \Omega_{\pm}^1}\bar w_{\pm}(\cdot,-1)^2d\bar \nu^{-1}\le
(1+Cr^2)\int_{\tilde \Omega_{\pm}^1}\tilde w_{\pm}^2d\nu_0^{-1}
\end{equation}
 and
\begin{equation}\label{55e2}
\int_{\bar \Omega_{\pm}^1}|\nabla_{\bar g}\bar
w_{\pm}(\cdot,-1)|^2d\bar \nu^{-1}\ge (1-Cr^2)\int_{\tilde
\Omega_{\pm}^1}|\nabla \tilde w_{\pm}|^2d\nu_0^{-1}.
\end{equation}
Recall that $d\nu_0^{-1}=\frac{1}{(4\pi)^{\frac
n2}}e^{-\frac{|x|^2}{4}}dx$.

 Beckner-Kenig-Pipher inequality (a proof of which can be found
in \cite{ck}) gives
\begin{equation}\label{55e3}
\frac{\int_{\tilde \Omega_+}|\nabla \tilde
w_+|^2d\nu_0^{-1}}{\int_{\tilde \Omega_+}\tilde w_+^2
d\nu_0^{-1}}+\frac{\int_{\tilde \Omega_-}|\nabla \tilde
w_-|^2d\nu_0^{-1}}{\int_{\tilde \Omega_-}\tilde w_-^2 d\nu_0^{-1}}
\ge 1.
\end{equation}
Therefore Proposition \ref{09428p1} follows immediately from
(\ref{55e1}), (\ref{55e2}) and (\ref{55e3}). $\Box$

\section{Proof of Proposition \ref{keyprop1} and Proposition \ref{keyprop2}}

First we observe that Proposition \ref{keyprop2} is a direct
consequence of Lemma \ref{harp1} and Lemma \ref{harp2} as in
\cite{cjk,edquist}. To prove Proposition \ref{keyprop1} we let

$$\tilde v_{\pm}(y,s_1)=\frac{1}{4^{2k}}w_{\pm}(4^ky,4^{2k}s_1),\quad (y,s_1)\in
\tilde \Omega_{\pm}. $$  It is easy to see that $\tilde \Omega_+$
and $\tilde \Omega_-$ are disjoint subsets of $B(0,\delta 4^k)$
where $\delta>0$ is small. Let $\tilde g_{ij}$ be the scaled metric,
$d\tilde \nu$ and $d\tilde \nu^s$ be the new measures. Let
$$\tilde \phi(r)=\frac 1{r^4}\int\!\!\!\int_{S_r}|\nabla_{\tilde
g}\tilde v_+|^2d\tilde \nu\int\!\!\!\int_{S_r}|\nabla_{\tilde
g}\tilde v_-|^2d\tilde \nu.$$ Also we set
$$\tilde A_{\pm}(r)=\int\!\!\!\int_{S_r}|\nabla_{\tilde
g}\tilde v_{\pm}|^2d\tilde \nu, \quad  \tilde
B_{\pm}(r)=\int_{\mathbb R^n}|\nabla_{\tilde g}\tilde v_{\pm}|^2
\tilde U(y,-r^2)dV_{\tilde g}. $$

We want to show that for $\frac 14\le r\le 1$, if $\tilde
A_{\pm}=\tilde A_{\pm}(1)$ are both large, then
\begin{equation}\label{429e1}
\tilde \phi'(r)\ge -C\tilde \phi(r)(\frac{1}{\sqrt{\tilde
A_+}}+\frac{1}{\sqrt{\tilde A_-}}+4^{-2k})
\end{equation}
where $C$ is independent of $k$. Once (\ref{429e1}) is established,
the integration of (\ref{429e1}) gives
\begin{equation}\label{429e2}
\tilde \phi(\frac 14)\le \tilde \phi(1)(1+C\delta_k),\quad
\delta_k=\frac{1}{\sqrt{\tilde A_+}}+\frac{1}{\sqrt{\tilde
A_-}}+4^{-2k}.
\end{equation}
Then by scaling, (\ref{429e2}) is equivalent to Proposition
\ref{keyprop1}.

So we are left with the proof of (\ref{429e1}).  By scaling, the
case for $\frac 14\le r\le 1$ can be treated as $r=1$. Then the
proof is very similar to the standard one:
$$\tilde \phi'(1)=-4\tilde A_+\tilde A_-+2\tilde B_+\tilde A_-+2\tilde
A_+\tilde B_-.$$

If $\tilde B_+\ge 2\tilde A_+$ or $\tilde B_-\ge 2\tilde A_-$,
$\tilde \phi'(1)\ge 0$. So we only assume $\tilde B_{\pm}\le 4\tilde
A_{\pm}$.

Now we apply (\ref{322e1}) to $\tilde v_{\pm}$ to get (using
$d\tilde \nu^{-1}=\tilde U(\cdot, 1)dV_{\tilde g}$)
\begin{eqnarray}
&&\tilde A_{\pm}\le C+C(\int_{\mathbb R^n}\tilde v_{\pm}^2d\tilde
\nu^{-1})^{\frac 12}+\frac 12\int_{\mathbb R^n}\tilde
v_{\pm}^2d\tilde
\nu^{-1}\nonumber\\
\label{429e3}
 &\le &C+\frac{C}{\sqrt{\lambda_{\pm}}}\sqrt{\tilde
B_{\pm}}+\frac 1{2\lambda_{\pm}}\tilde B_{\pm}.
\end{eqnarray}
Note that we can assume $\lambda_+$ and $\lambda_-$ are both
positive, because if, say $\tilde v_+\equiv 0$, we obtain from the
first line of (\ref{429e3}) that
$$\tilde A_{+}\le C,$$
which is a contradiction to the largeness of $\tilde A_+$.

From (\ref{429e3}) we see that if $\lambda_+\ge 2$ or $\lambda_-\ge
2$, (\ref{429e1}) is established easily. Therefore we assume
$\lambda_{\pm}\le 2$. In this case we obtain from (\ref{429e3}) that
\begin{equation}\label{429e4}
2\lambda_+\tilde A_+\le C+\frac{C}{\sqrt{\lambda_+}}\sqrt{\tilde
B_+}+\tilde B_+.
\end{equation}
There is a similar equation for $2\lambda_-\tilde A_-$. Multiplying
$\tilde A_-$ to (\ref{429e4}), $\tilde A_+$ to the corresponding
equation, we have, by adding these two equations
\begin{eqnarray}\label{429e5}
2(\lambda_++\lambda_-)\tilde A_+\tilde A_-&\le & C(1+\tilde
A_-\sqrt{\frac{\tilde B_+}{\lambda_+}}+\tilde A_+\sqrt{\frac{\tilde
B_-}{\lambda_-}})\\
&&+\tilde B_+\tilde A_-+\tilde B_-\tilde A_+.\nonumber
\end{eqnarray}
Proposition \ref{09428p1} gives $\lambda_++\lambda_-\ge 1-C
4^{-2k}$. Using this in (\ref{429e5}) we obtain (\ref{429e1}).
Proposition \ref{keyprop1} is established. $\Box$

Theorem \ref{thm1} follows from Proposition \ref{keyprop1} and
Proposition \ref{keyprop2} as in \cite{cjk} and \cite{edquist}.

\bigskip

If further information on the growth of $u_{\pm}$ is known near the
origin, then the behavior of $\phi(r)$ can be made more precise.
This is the observation in \cite{cjk,edquist}.

\begin{thm}\label{thm2}
Let $u_{\pm}$,$w_{\pm}$,$\chi$ be the same as in Theorem \ref{thm1},
suppose in addition that
$$|u_{\pm}(x,s)|\le C_{\epsilon}(|x|^2+|s|)^{\frac{\epsilon}2}$$
for $(x,s)\in Q_{\delta_p}$ and $\epsilon\in (0,1]$. Then
$$\phi(r)\le (1+\rho^{\epsilon})\phi(\rho)+C_M\rho^{\epsilon},\quad
0<r\le \rho\le \delta_p/4$$ where $C_M$ depends on $n,\Lambda$,
$M=\|u_+\|_{L^2(Q_{\delta_p})}+\|u_-\|_{L^2(Q_{\delta_p})}$,$\chi$,$\epsilon$.
\end{thm}

Since the proof of Theorem \ref{thm2} is similar to its analogue
in \cite{edquist}, we leave the detail to the interested readers.
It is worthwhile to point out here that  the only difference in
the proof comes from the correction term in the
Beckner-Kenig-Pipher inequality, \cite{beckner}. The readers can
easily see that at this point the extra term does not cause
further difficulties in the argument.

\end{document}